\documentclass[11pt,a4paper]{article} 

\usepackage{amsmath,amssymb,graphicx}
\newtheorem{theoreme}{Theorem}[section]

\newtheorem{define}[theoreme]{Definition}
\newtheorem{propal}[theoreme]{Proposition}

\newenvironment{proof2}
    {\begin{list}{{\bf Proof}}
                 {\setlength{\labelwidth}{0cm}       %{3cm}
                  \setlength{\leftmargin}{0em}}    %{3.6em}}
                 \item}
    {$\blacksquare$\end{list}}
\newenvironment{remarque}{\par \addvspace{\bigskipamount} 
                               \refstepcounter{theoreme}
                               {\bf Remark \thetheoreme :~}}
                         {\par \addvspace{\bigskipamount}}
\def \disp {\displaystyle}

\def \R {{\mathbb R}}
\def \C {{\cal C}}
\def \eps {{\varepsilon}}

\def \D {{\cal D}}

\def \E {{{\cal E}}}

\newcommand{\Const}{{\cal F}} 
\newcommand{\Conv}{{\textnormal{Conv}}}
\newcommand{\dist}{{\delta}}

\def \un {{\ensuremath{\textrm{\textnormal{1}}\hspace{-0.7ex}
                       \textrm{\textnormal{I}}}}}

\def \bx {{\bf x}} \def \nx {{\breve{x}}}
\def \by {{\bf y}} 
\def \bz {{\bf z}}  
\def \bv {{\bf v}} \def \nv {{\breve{v}}}
\def \bX {{\bf X}} \def \nX {{\breve{X}}}
\def \bW {{\bf W}} \def \nW {{\breve{W}}}
  
\def \bl {{\bf l}}
\def \bL {{\bf L}} 
\def \bsigma {{{\boldsymbol\sigma}}} \def \nsigma {{\breve{\sigma}}}
\def \bb {{\bf b}}  

\def \rm {r_-}
\def \rp {r_+}
\newcommand \grav {{\varphi}}

\newcommand{\obl}{{\boldsymbol\theta}}

\author{\Large
Myriam \textsc{Fradon\footnote{U.F.R. de Math\'ematiques, CNRS U.M.R. 8524, %\newline \hspace*{1.4em}
                               Universit\'e Lille1, 59655 Villeneuve d'Ascq Cedex, France \newline \hspace*{1.4em}
                               Myriam.Fradon@univ-lille1.fr}
              }
       }
\date{}
\title{\bf\LARGE SDE with oblique reflection on domains defined by multiple constraints}

\begin{document}%%%%%%%%%%%%%%%%%%%%%%%%%%%%%%%%%%%%%%%%%%%%%%%%%%%%%%%%%%%%%%%%%

\maketitle
\centerline{\textbf{Abstract}}
We present simple assumptions on the constraints defining a hard core dynamics for the associated reflected stochastic differential equation to have a unique strong solution. Time-reversibility is proven for gradient systems with normal reflection, or oblique reflection with a fixed oblicity matrix. An application is given concerning the clustering at equilibrium of particles around a large attractive sphere.
\vspace{3mm}

\noindent
{\bf AMS 2000 subject classification:} 60K35, 60J55, 60H10.\\
% 60H10 	Stochastic ordinary differential equations
% 60J55    	Local time and additive functionals
% 60J60    	Diffusion processes
% 60K35    	Interacting random processes; statistical mechanics type models; percolation theory
{\bf Keywords:} Stochastic Differential Equation, hard core interaction, reversible measure, local time.\\%normal reflection

%>>>>>>>>>>>>>>>>>>>>>>>>>>>>>>>>>>>>>>>>>>>>>>>>>>>>>>>>>>>>>>>>>>>>>>>>>>>>>>>>>
\section{Introduction}
%>>>>>>>>>>>>>>>>>>>>>>>>>>>>>>>>>>>>>>>>>>>>>>>>>>>>>>>>>>>>>>>>>>>>>>>>>>>>>>>>>
%>>>>>>>>>>>>>>>>>>>>>>>>>>>>>>>>>>>>>>>>>>>>>>>>>>>>>>>>>>>>>>>>>>>>>>>>>>>>>>>>>

Since the first works of Skorokhod \cite{Skorokhod1et2} on existence and uniqueness for strong solutions of reflected stochastic differential equations, many authors have investigated this type of equation and extended his results on half-spaces to more general domains: convex sets (Tanaka \cite{Tanaka}), admissible sets (Lions-Sznitman \cite{LionsSznitman}), domains satisfying only the Uniform Exterior Sphere and the Uniform Normal Cone conditions (Saisho \cite{SaishoSolEDS}) or some weaker version of these conditions (Dupuis and Ishii \cite{DupuisIshii}). The question of equilibrium states of the reflected process (construction of time-reversible initial measures) has also been investigated (see e.g. \cite{SaishoTanakaSymmetry}).

All these existence, uniqueness and reversibility results were obtained under some smoothness assumptions on the boundary of the domain. Typically, the existence of at least one normal inward vector at each point of the boundary is a necessary condition to define the normal reflection direction.

In most cases, the domain in which the process has to live is defined by constraints which are physically natural rather than by its geometrical properties as a subset of some Euclidean space. For example, consider a system of $n$ identical hard spheres with radius $r$ in $\R^d$. The domain in which they evolve is the set of configurations $(x_i)_{1 \le i \le n}$ satisfying the constraints $|x_i-x_j|>2r$ (i.e. the distance between the centers of any two spheres is larger than twice their radius); the geometrical description is much more complicated: the complementary set in $\R^{nd}$ of some star-convex subset whose boundary can be locally approximated by a tangent sphere and a cone.

Unfortunately, for reflected processes in dimension larger than three, the geometrical properties of the domain are not that obvious from the physical constraints. In the already mentioned $nd$-dimensional example of a finite system of hard spheres, the main part of the paper \cite{SaishoTanakaBrownianBalls} is devoted to proving that the set of allowed ball configurations satisfies the Uniform Exterior Sphere and Uniform Interior Cone property. In \cite{SaishoTanemuraBrownianParticles} and \cite{FradonBrownianGlobules} too, a meticulous and extensive geometrical study has to be performed before the stochastic analysis of the dynamics. 

We present in this note a constraint-based assumption for existence, uniqueness and reversibility results for Skorokhod equations.

Our aim is to deal with assumptions as simple and physically natural as possible, even if they are not the weakest ones.
Moreover, we treat the delicate case of oblique reflection.
It will be the basis of a general approach for analysing stochastic dynamics under constraints.
We also prove that for well-chosen diffusion coefficients the oblique reflected gradient diffusion admits a reversible measure.

This note is divided in two parts. 

The first part (section \ref{sect_Skorohod_eq}) exhibits a new compatibility criterion for constraints. If it is satisfied, then the reflected stochastic differential equation admits a unique strong solution, and this solution is time-reversible in the case of a gradient system whose reflection direction is consistent with its diffusion coefficient. 

In section \ref{sect_examples} we present an application.
We consider the behaviour of many particles around a sphere called planet.
%As in \cite{FradonBrownianGlobules}, 
The particles are spherical with random radii oscillating between a minimum and a maximum value. Their motion is perturbated by collisions into other particles and into the planet. The planet also generates a gravitational field which has a smooth attractive influence on the particles. 
We are interested in the equilibrium situation:
Are the particles scattered in a wide area as a typical configuration, or do they tend to cluster around the planet~? We prove that at equilibrium and for low temperature all particles are as close as possible to the planet, all located beneath some altitude, with high probability.

%>>>>>>>>>>>>>>>>>>>>>>>>>>>>>>>>>>>>>>>>>>>>>>>>>>>>>>>>>>>>>>>>>>>>>>>>>>>>>>>>>
\section{Reflected stochastic differential equation under multiple constraints}
\label{sect_Skorohod_eq}
%>>>>>>>>>>>>>>>>>>>>>>>>>>>>>>>>>>>>>>>>>>>>>>>>>>>>>>>>>>>>>>>>>>>>>>>>>>>>>>>>>
%>>>>>>>>>>>>>>>>>>>>>>>>>>>>>>>>>>>>>>>>>>>>>>>>>>>>>>>>>>>>>>>>>>>>>>>>>>>>>>>>>

We are interested in a process living in the closure of a domain $\D$. This domain is defined by a finite set $\Const$ of smooth $\R$-valued constraint functions on $\R^d$:
$$
\D=\left\{ \bx \in \R^d;~f(\bx)>0 \text{ for each } f\in\Const \right\}.
$$
$\D$ is an intersection of smooth sets (arbitrary many of them, provided they are in finite number), so its boundary is a finite union of smooth boundaries:
$$
\partial\D=\bigcup_{f\in\Const} \left\{ \bx \in \D;~f(\bx)=0 \right\}.
$$
Since we want the process to be reflected on the boundary of  $\D$, we have to assume some regularity on the functions in $\Const$.
The reflection at any point $\bx\in\partial\D$ occurs either in the inward normal direction $\nabla f(\bx)$ or in an oblique direction depending only on the normal direction. So we have to suppose the existence of a direction which is normal to the boundary:
$\nabla f(\bx) \neq 0$ for each $\bx\in\overline{\D}$ such that $f(\bx)=0$.
We actually assume something more: the first derivatives of the functions of $\Const$ should admit some positive uniform lower bound, their second derivatives should be uniformly bounded and, most important, we suppose that the boundary of each single-constraint set $\left\{ \bx~;~ f(\bx) > 0 \right\}$ crosses the boundaries of the other single-constraint sets at "not too sharp an angle". To be more precise, we have to exclude infinitely sharp "thorns" which contain a point admiting inward normal vectors with opposite directions. This is what we call \emph{compatibility} between the constraints:
\begin{define}\label{def_compatibilite}
Let $\Const$ be a finite set of $\R$-valued $\C^2$-functions on $\R^d$.
These functions are called \emph{compatible constraints} if 
\begin{itemize} 
\item $\disp \D:=\left\{ \bx \in \R^d;~f(\bx)>0 \text{ for each } f\in\Const \right\}$ is a non-empty connected set~;
\item for each $f\in\Const$,~~$\disp \inf\{ |\nabla f(\bx)|;~\bx\in\overline{\D},~ f(\bx)=0 \} >0$ ~and~
      $\disp \sup\{ |D^2 f(\bx)|;~\bx\in\R^d \} < +\infty $~;
\item $\disp \inf_{\bx\in\partial\D} \dist(0,\Conv(\bx)) > 0$ \\
where $\Conv(\bx)$ is the convex hull of the unit normal vectors to the boundaries at point $\bx$:
$$
\Conv(\bx)=\left\{ \sum_{f\in\Const,f(\bx)=0} c_f \frac{\nabla f(\bx)}{|\nabla f(\bx)|} 
                               \text{ ~~s.t.~~ } c_f \ge 0 \text{ and} \sum_{f\in\Const,f(\bx)=0} c_f =1  \right\}.
$$
\end{itemize} 
\end{define}
Here and in the sequel, $\dist$ denotes the Euclidean distance in $\R^d$, $|\by|$ denotes the Euclidean norm of vector $\by$, and  $|M|=\sup\{ |M\by|/|\by|~;~\by\in\R^d \}$ denotes the norm of the matrix $M$. Lebesgue measure is denoted by $d\bx$. 

The next main theorem states that our compatibility definition provides a convenient assumption to ensure the existence of a reflected process within a set defined by constraints. In most models, for the sake of simplicity, the reflection direction is the inward \emph{normal} direction on the boundary.
Here we consider a general \emph{oblique} reflection, as in the case treated in \cite{FradonRoellyGlobulesInfinis} or in Section \ref{sect_examples}. We state the result with a fixed deviation $\obl~^t\obl$ from the normal direction. $^t\obl$ denotes the transposed matrix. Normal reflection corresponds to the special case $\obl=I_d$. 

\begin{theoreme}{\bf (Existence and uniqueness of the reflected process)}\label{thPrincipal} \\
Let $\obl$ be a fixed $d \times d$ invertible matrix and $\Const$ be a set of compatible constraints, with $\D=\bigcap_{f\in\Const}\left\{ \bx \in \R^d;~f(\bx)>0 \right\}$ the corresponding subset of $\R^d$.
If $\bsigma:\overline{\D}\longrightarrow\R^{d^2}$ and $\bb:\overline{\D}\longrightarrow\R^d$ are bounded Lipschitz continuous functions on $\disp \overline{\D}$, then the reflected stochastic differential equation:
\begin{equation}\label{eqPrincipale}
\bX(t)=\bx +\int_0^t \bsigma(\bX(s)) d\bW(s) +\int_0^t \bb(\bX(s)) ds +\sum_{f\in\Const} \int_0^t \obl~^t\obl \nabla f(\bX(s)) dL_f(s) 
\end{equation}
has for each starting point $\bx\in\D$ a unique strong solution in $\disp \overline{\D}$, where the local times $L_f$ satisfy
$\disp L_f(\cdot) = \int_0^\cdot \un_{f(\bX(s))=0}~dL_f(s)$.
\end{theoreme}
In this theorem "Strong uniqueness of the solution" stands for strong uniqueness, in the sense of \cite{IkedaWatanabe} chap.IV def.1.6, of the process $\bX$, not of the local times $L_f$.

If $\bsigma$ is constant and the drift $\bb$ is a gradient, then the equation also admits a time-reversible measure $\mu$ (i.e. the distribution of the solution with initial measure $\mu$ is invariant under the transformation $\big(\bX(\cdot),(L_f(\cdot))_{f\in\Const}\big)\longrightarrow\big(\bX(T-\cdot),(L_f(T-\cdot)-L_f(T))_{f\in\Const}\big)$ for each $T>0$):

\begin{theoreme}{\bf (Reversibility of the reflected gradient process)}\label{thReversibilite} \\
Let $\obl$ be a fixed $d \times d$ invertible matrix and $\Const$ a set of compatible constraints. If $\Phi$ is a $\C^2$-function on $\R^d$ with bounded derivatives, then the solution of
\begin{equation}\label{eqObliqueReflechie}
\bX(t)=\bX(0) +\obl\bW(t) -\frac{1}{2} \int_0^t \obl~^t\obl \nabla\Phi(\bX(s)) ds +\sum_{f\in\Const} \int_0^t \obl~^t\obl \nabla f(\bX(t)) d\bL_f(s) 
\end{equation}
admits $d\mu(\bx)=\un_{\D}(\bx) e^{-\Phi(\bx)} d\bx$ as a time-reversible measure.
\end{theoreme}

\begin{remarque}\label{remEquivalence}
In %the compatibility 
definition \ref{def_compatibilite}, the condition $\disp \inf_{\bx\in\partial\D} \dist(0,\textnormal{Conv}) > 0$ is equivalent to:
$$
\exists \beta_0 >0 \quad \forall\bx\in\partial\D \quad \exists \bv\neq 0, \quad
\forall f\in\Const \text{ s.t. } f(\bx)=0 \quad \bv.\nabla f(\bx) \ge \beta_0 |\bv|~|\nabla f(\bx)|
$$
where the dot denotes the Euclidean scalar product.

Though this statement is longer and apparently more difficult to obtain than an uniform lower bound on the norms of the convex combinations, it is in some sense more intuitive. It states the existence of cones (with vertex $\bx$, axis $\bv$ and aperture $2\arccos \beta_0$ ) which contain all the inward normal vectors given by the constraints at point $\bx$. The positivity condition ensures that these cones do not degenerate into half-spaces.
This condition is easier to check in some concrete situations (see section \ref{sect_examples}).
\end{remarque}
\begin{remarque}{\bf (Stability of the compatibility property)}\label{remStabilite}\\
Let $\Const$ be a set of compatible constraints on $\R^d$.
\begin{itemize}
\item
If $\obl$ is a $d \times d$ invertible matrix, the transformed constraints $\{ f(\obl~\cdot);~f\in\Const \}$ are compatible.
\item
If all constraints disregard one of the $d$ coordinates, then $\Const$ induces a set of compatible constraints on $\R^{d-1}$; that is if $f(x_1,\cdots,x_{d-1},x_d)=f(x_1,\cdots,x_{d-1},0)$ for each $f$ in $\Const$ and each $\bx=(x_1,\cdots,x_d)$ in $\R^d$ then 
$\disp \left\{ \underline{f}:{\small \begin{array}{l}\R^{d-1} \longrightarrow \R                   \\
                                             \underline\bx \longmapsto f(\underline\bx,0)
                                     \end{array} } ;~f\in\Const \right\}$
is compatible.
\end{itemize}
\end{remarque}
The remaining of this section is devoted to the proofs of the above results. We first prove remarks \ref{remEquivalence} and \ref{remStabilite} which will be useful in the other proofs, and then proceed to theorems \ref{thPrincipal} and \ref{thReversibilite}.

\begin{proof2} {\bf of remark \ref{remEquivalence}} \\
The third compatibility condition is $\disp \exists \beta_0 >0 ~~ \forall\bx\in\partial\D ~~ \dist(0,\Conv(\bx)) \ge \beta_0$. 
The condition in remark \ref{remEquivalence} can be rewritten as
$\disp \exists \beta_0 >0 ~ \forall\bx\in\partial\D ~ \max_{\bv\neq 0} \min \left\{ \frac{\bv}{|\bv|}.\frac{\nabla f(\bx)}{|\nabla f(\bx)|};~ f\in\Const,~ f(\bx)=0 \right\} \ge \beta_0$.
Thus it suffices to prove that for each $\bx\in\partial\D$
$$
\dist(0,\Conv(\bx))=\max_{|\bv|=1} \min \left\{ \bv.\frac{\nabla f(\bx)}{|\nabla f(\bx)|};~ f\in\Const,~ f(\bx)=0 \right\}.
$$
The lower bound on $\dist(0,\Conv(\bx))$ follows from the inequality 
$|\by| \ge \by.\bv \ge \min_{f\in\Const,f(\bx)=0} \frac{\nabla f(\bx)}{|\nabla f(\bx)|}.\bv$,
which holds for every unit vector $\bv$ and every $\by\in\Conv(\bx)$ because families $(c_f)$ of non-negative numbers summing up to $1$ satisfy:
$$
\left( \sum_{f,f(\bx)=0} c_f \frac{\nabla f(\bx)}{|\nabla f(\bx)|} \right).\bv 
\ge \left( \sum_{f,f(\bx)=0} c_f \right) \min_{f,f(\bx)=0} \frac{\nabla f(\bx)}{|\nabla f(\bx)|}.\bv
$$
Since the convex hull $\Conv(\bx)$ is a closed set, it contains an element $\bz$ with minimal norm: $|\bz|=\dist(0,\Conv(\bx))$.
For each $f$ satisfying $f(\bx)=0$ and for each positive $\eps$, the convex combination
$\frac{1}{1+\eps} \left( \bz + \eps~\frac{\nabla f(\bx)}{|\nabla f(\bx)|} \right)$
belongs to the convex hull, hence its norm can not be smaller than $|\bz|$:
$$
|\bz|^2 + \eps^2 + 2~\eps~\bz.\frac{\nabla f(\bx)}{|\nabla f(\bx)|} \ge (1+\eps)^2~|\bz|^2
\qquad\text{ i.e. }\qquad
\eps + 2~ \bz.\frac{\nabla f(\bx)}{|\nabla f(\bx)|} \ge (2 + \eps) |\bz|^2
$$
This proves that
$\disp \frac{\bz}{|\bz|}.\frac{\nabla f(\bx)}{|\nabla f(\bx)|} \ge |\bz|=\dist(0,\Conv(\bx))$ 
and provides the upper bound on $\dist(0,\Conv(\bx))$.
%for each $f$ and some $\bz$, and consequently $\disp d(0,\Conv(\bx)) \le \max_{|\bv|=1} \min_{f,f(\bx)=0} \bv.\frac{\nabla f(\bx)}{|\nabla f(\bx)|}$.
\end{proof2}

% on utilise les formules de derivation des fonctions de plusieurs variables (verifiees !) :
% \nabla f= ^t Df  qui traduit  Df(x)(y) = (\nabla f(x)).y  (Df(x) vecteur ligne, \nabla f(x) vecteur colonne)
% D(f(\obl~\cdot))=(Df)(\obl~\cdot)\obl
% \nabla(f(\obl~\cdot))= ^t\obl (\nabla f|)(\obl~\cdot)
% D^2(f(\obl~\cdot))=^t\obl(D^2 f)(\obl~\cdot)\obl

\begin{proof2} {\bf of remark \ref{remStabilite}} \\
Let us prove the compatibility of the set $\Const^\obl=\{ g(\cdot)=f(\obl~\cdot);~f\in\Const \}$ of transformed constraints.
Since matrix $\obl$ is invertible, $\obl^{-1}\D=\left\{ \by \in \R^d;~ \forall f\in\Const~f(\obl \by) > 0 \right\}$ is a non-empty connected set as continuous image of the non-empty connected set $\D$.
$\obl$ also transforms the bounds on the $f$'s into bounds on the $g$'s.
%$$\begin{array}{l} \disp
%\inf_{\{g=0\}\cap\obl^{-1}\D}{\nabla g} =\inf\{ |^t\obl \nabla f(\bx)|;~ \bx\in\D,~ f(\bx)=0 \}
%\ge \frac{1}{|^t\obl^{-1}|} \inf\{ |\nabla f(\bx)|;~ \bx\in\D,~ f(\bx)=0 \} >0 \\ \disp
%||D^2 g||_\infty =\sup\{ |^t\obl (D^2 f)(\obl \by) \obl|;~\by\in\R^d \} \le |^t\obl| \sup\{ |D^2 f(\bx)|;~\bx\in\R^d \} |\obl|< +\infty 
%\end{array}
%$$
Remark \ref{remEquivalence} with $\bv$ replaced by $\obl \bv$ provides the existence of some positive $\beta_0$ such that:
$$
\forall\bx\in\partial\D \quad \exists \bv\neq 0, \quad
\forall f\in\Const \text{ s.t. } f(\bx)=0 \quad \bv.^t\obl \nabla f(\bx) \ge \beta_0 |\obl \bv|~|\nabla f(\bx)|.
$$
Replacing $\bx$ by $\obl\by$ we obtain:
$$
\forall\by\in\partial(\obl^{-1}\D) \quad \exists \bv\neq 0, \quad
\forall g\in\Const^\obl \text{ s.t. } g(\by)=0 \quad 
\bv.\nabla g(\by) \ge \beta_0 |\obl \bv|~|^t\obl^{-1} \nabla g(\by)| 
                  \ge \beta_0 \frac{|\bv|}{|\obl^{-1} |} \frac{|\nabla g(\by)|}{|^t\obl|}
$$
Thanks to remark \ref{remEquivalence} with $\beta_0'=\frac{\beta_0}{|\obl^{-1}|~|^t\obl| }$, this proves that 
$\Const^\obl$ is a set of compatible constraints.

In order to prove the second part of remark \ref{remStabilite}, we now assume that $f(\underline{\bx},x_d)=f(\underline{\bx},0)$ for each $f$ in $\Const$ and each $(\underline{\bx},x_d)$ in $\R^d$. 
%We want to prove that the $\underline{f}$'s defined on $\R^{d-1}$ by $\underline{f}:\underline{\bx} \longmapsto f(\underline{\bx},0)$ are compatible on $\R^{d-1}$.
The set $\D=\{ \bx\in\R^d;~f(\bx)>0 \}$ is equal to $\underline{\D} \times \R$ where $\underline{\D}=\{\bz\in\R^{d-1};~\underline{f}(\bz)>0\}$ is a non empty connected set as a projection of a non-empty connected set.
The lower bound on $\nabla f$ and the upper bound on $D^2 f$ transfer to $\underline{f}$ because $\nabla f=(\nabla\underline{f},0)$ and $|D^2 f(\bx)|=|D^2 f(x_1,\cdots,x_{d-1},0)|$.
%The $f$'s in $\Const$ are constant as functions of the last coordinate in $\R^d$ thus $\nabla f=(\nabla\underline{f},0)$ and
%$$
%\inf\{ |\nabla\underline{f}(\bz)|;~ \bz\in\underline{\D},~ \underline{f}(\bz)=0 \ =\inf\{ |\nabla f(\bz,0)|;~ (\bz,0)\in\D,~ f(\bz,0)=0 \} >0
%$$
%For the same reason, $|D^2 f(\bx)|=|D^2 f(x_1,\cdots,x_{d-1},0)|$ and
%$$
%\sup\{ |D^2 \underline{f}(\bz)|;~ \bz\in\R^{d-1} \} = \sup\{ |D^2 f(\bx)|;~ \bx\in\R^d \} < +\infty 
%$$
From the compatibility of $\Const$, we also get the existence of a positive $\beta_0$ such that for each $\bx\in\partial\D$ there exists a unit vector $\bv$ satisfying $\bv.\nabla f(\bx) \ge \beta_0 |\nabla f(\bx)|$ for each function $f\in\Const$ vanishing at point $\bx$.
%$$
%\exists \beta_0 >0 \quad \forall\bx\in\partial\D \quad \exists \bv\neq 0, \quad
%\forall f\in\Const \text{ s.t. } f(\bx)=0 \quad \bv.\nabla f(\bx) \ge \beta_0 |\bv|~|\nabla f(\bx)|
%$$
%$\bv.\nabla f(\bx)$ is positive and 
The last coordinate of $\nabla f(\bx)$ vanishes, hence $\underline{\bv}=(v_1,\cdots,v_{d-1}) \neq 0$. Since $\partial\D = \partial\underline{\D} \times \R$ we obtain the compatibility of the $\underline{f}$'s:
$$
\exists \beta_0 >0 \quad \forall\bz\in\partial\underline{\D} \quad \exists \underline{\bv}\neq 0, \quad
\forall f\in\Const \text{ s.t. } \underline{f}(\bz)=0 \quad \underline{\bv}.\nabla\underline{f}(\bz) \ge \beta_0 |\underline{\bv}|~|\nabla\underline{f}(\bz)|
$$
\end{proof2}

\begin{proof2} {\bf of theorems \ref{thPrincipal} and \ref{thReversibilite}}\\
{\bf The case of normal reflection:} we assume here that $\obl=I_d$.
According to corollary 3.6 of \cite{FradonBrownianGlobules}, equation (\ref{eqPrincipale}) has a unique strong solution as soon as $\D$ satisfies the four assumptions of the inheritance criterion for Uniform Exterior Sphere and Uniform Normal Cone conditions (proposition 3.4 in \cite{FradonBrownianGlobules}). We will check these four assumptions in the unusual order (i) (ii) (iv) (iii)  because some parameter appearing in (iii) depends on a parameter defined in (iv). We use the notations:
$\disp \underline{\nabla f}:=\inf\{ |\nabla f|(\bx);~ \bx\in\overline{\D},~ f(\bx)=0 \}$ and 
$\disp || D^2f ||_\infty := \sup\{ |D^2 f|(\bx);~ \bx\in\R^d \}$.

{\it Assumption (i):}
We have to prove that $\left\{ \bx \in \R^d~;~ f(\bx) \ge 0 \right\}$ has $\C^2$ boundary in $\overline{\D}$ for each constraint $f$. 
Let us fix $\bx\in\overline{\D}$ such that $f(\bx)=0$. By definition of the constraint functions, $\nabla f(\bx) \neq 0$, that is we can choose an index $k$ such that $\nabla_k f(\bx) \neq 0$. For simplicity sake, we assume that $\nabla_d f(\bx) > 0$ (the idea easily adapts to $k \neq d$ and to negative partial derivatives). 
%Applying the implicit function theorem to the $\C^2$-function $f$, we obtain the existence of a neighborhood $V$ of $(x_1,\ldots,x_{d-1})$, a neighborhood $U'$ of $x_d$, an $\eps>0$ and a $\C^2$-function $h$ such that $f(y_1,\ldots,y_d)=z_d \Leftrightarrow y_d=h(y_1,\ldots,y_{d-1},z_d)$ for $(y_1,\ldots,y_{d-1}) \in V$, $y_d \in U'$ and $|z_d|<\eps$. Since $\nabla_d f(\bx)$ is positive, $U'$ can be chosen small enough for $f$ to be increasing on $U'$ as a function of $y_d$. Then $h$ is increasing as a function of $z_d$:
%$f(y_1,\ldots,y_d)=0 \Leftrightarrow y_d=h(y_1,\ldots,y_{d-1},0)$ and 
%$f(y_1,\ldots,y_d)>0 \Leftrightarrow y_d>h(y_1,\ldots,y_{d-1},0)$. 
% Thus
Applying the implicit function theorem to the $\C^2$-function $f$, we obtain the existence of a neighborhood $V$ of $(x_1,\ldots,x_{d-1})$, a neighborhood $U'$ of $x_d$ and an increasing $\C^2$-function $h$ such that the $\C^2$-diffeomorphism $(y_1,\ldots,y_d) \longmapsto (y_1,\ldots,y_{d-1},f(y_1,\ldots,y_d))$ maps 
$\left\{ \by \in V \times U',~ f(\by) \ge 0 \right\}$ to 
$\left\{ (y_1,\ldots,y_{d-1},z_d) \in V \times U';~ z_d \ge h(y_1,\ldots,y_{d-1},0) \right\}$.
Hence, the subset $\left\{ \bx \in \R^d;~ f(\bx) \ge 0 \right\}$ has $\C^2$ boundary in $\overline{\D}$, and its inward normal direction at point $\bx$ is $\disp \frac{\nabla f(\bx)}{|\nabla f(\bx)|}$.

{\it Assumption (ii):} Let us prove that $\left\{ \bx \in \R^d;~ f(\bx) \ge 0 \right\}$ satisfies the Uniform Exterior Sphere property restricted to $\overline{\D}$:
according to definition 3.1 in \cite{FradonBrownianGlobules}, we have to prove that there exists some positive $\alpha_f$ such that, for each $\bx\in\overline{\D}$ satisfying $f(\bx)=0$, one has
\begin{equation}\label{eqUES}
\forall \by \text{ s.t. } f(\by) \ge 0 \qquad (\by-\bx).\frac{\nabla f(\bx)}{|\nabla f(\bx)|}+\frac{1}{2\alpha_f}|\by-\bx|^2 \ge 0.
\end{equation}
Let us fix $\bx\in\overline{\D}$ on which $f$ vanish. Taylor formula gives:
$$
\nabla f(\bx).(\by-\bx) + \frac{1}{2} (\by-\bx). D^2f(\bx+c^*(\by-\bx)) (\by-\bx) = f(\by) 
$$
for each $\by\in\R^d$, with some $c^* \in [0;1]$ depending on $\by$ and $\bx$.
In particular, for $\by$ such that $f(\by) \ge 0$ we obtain
$\disp (\by-\bx).\frac{\nabla f(\bx)}{|\nabla f(\bx)|} +\frac{||D^2f||_\infty}{2|\nabla f(\bx)|} |\by-\bx|^2 \ge 0$
which gives (\ref{eqUES}) with $\disp \alpha_f=\frac{\underline{\nabla f}}{||D^2f||_\infty}$.

{\it Assumption (iv):} %This is the main compatibility assumption. 
We have to prove the existence of some $\beta_0>0$ such that for each $\bx\in\partial\D$ there exists a unit vector $\bl_\bx^0$ satisfying $\bl_\bx^0.\nabla f(\bx) \ge \beta_0 |\nabla f|$ for each constraint such that $f(\bx)=0$. But this has already been done in remark \ref{remEquivalence} with $\beta_0=\inf_{\bx\in\partial\D} d(0,\Conv(\bx))$ and $\bl_\bx^0=\frac{\bz}{|\bz|}$ for some $\bz$ with minimal norm in $\Conv(\bx)$.

{\it Assumption (iii):} We have to prove that each set $\left\{ \bx \in \R^d;~ f(\bx) \ge 0 \right\}$ satisfies the Uniform Normal Cone property restricted to $\overline{\D}$ with constant $\beta_f$ smaller than $\beta_0^2/2$. Using Taylor formula for the derivative of $f$ yields to $\disp \nabla f(\by)=\nabla f(\bx)+ D^2f(\bx+c^*(\by-\bx)) (\by-\bx)$ for some $c^* \in [0;1]$ depending on $\by$ and $\bx$. 
%Taking the scalar product and dividing by the norms, 
We obtain for $\bx$ and $\by$ on which $f$ vanish:
$$
\frac{\nabla f(\bx).\nabla f(\by)}{|\nabla f(\bx)||\nabla f(\by)|}
=\frac{|\nabla f(\bx)|}{|\nabla f(\by)|} 
 +\frac{\nabla f(\bx).D^2f(\bx+c^*(\by-\bx)) (\by-\bx)}{|\nabla f(\bx)||\nabla f(\by)|}
$$
Since $|\nabla f(\bx)| \ge |\nabla f(\by)|-|D^2f(\bx+c^*(\by-\bx)) (\by-\bx)|$ the right hand side is not smaller than:
$$
1 - \frac{|D^2f(\bx+c^*(\by-\bx)) (\by-\bx)|}{|\nabla f(\by)|} 
  -\frac{|\nabla f(\bx)||D^2f(\bx+c^*(\by-\bx)) (\by-\bx)|}{|\nabla f(\bx)||\nabla f(\by)|}
$$
that is 
$\disp \frac{\nabla f(\bx)}{|\nabla f(\bx)|}.\frac{\nabla f(\by)}{|\nabla f(\by)|} \ge 1 - 2\frac{||D^2f||_\infty}{\underline{\nabla f}} |\by-\bx|$.
As a consequence, for any $\beta_f\in ]0,1[$ one can choose a $\delta_f>0$ small enough such that for each $\bx\in\overline{\D}$ satisfying $f(\bx)=0$ and each $\by\in\overline{\D}$ satisfying $f(\by)=0$ and $|\by-\bx|\le\delta_f$ one has
$\disp \frac{\nabla f(\by)}{|\nabla f(\by)|}.\frac{\nabla f(\bx)}{|\nabla f(\bx)|} \ge \sqrt{1-\beta_f^2}$.

This proves that $\left\{ \bx \in \R^d;~ f(\bx) \ge 0 \right\}$ satisfies the Uniform Normal Cone property restricted to $\overline{\D}$ with any constant $\beta_f\in]0;1[$, in particular with $\beta_f<\beta_0^2/2$ as requested.

To complete the proof of theorems \ref{thPrincipal} and \ref{thReversibilite} for $\obl=I_d$, we proceed as in the proof of theorem 3.3 in \cite{FradonBrownianGlobules}, replacing the probability measure $d\mu(\bx)=\frac{1}{Z} \un_{\D}(\bx) e^{-\Phi(\bx)} d\bx$ in that proof by the ($\sigma$-finite but maybe unbounded) measure $\mu$ defined by $d\mu(\bx)=\un_{\D}(\bx) e^{-\Phi(\bx)} d\bx$. Girsanov theorem yields the density of the distribution of the process with initial measure $\mu$ with respect to the distribution of reflected Brownian motion with Lebesgue measure as initial measure. Since both this density and the distribution of reflected Brownian motion starting from Lebesgue measure are time-reversal invariant, we obtain the reversibility of the solution with initial measure $\mu$. 

{\bf The case of oblique reflection:}
Let us check that the results obtained in the normal reflection case $\obl=I_d$ transfer to the case of any invertible matrix $\obl$.
Using the notation $\bX^\obl=\obl^{-1}\bX$, existence and uniqueness for equation (\ref{eqPrincipale}) is equivalent to  existence and uniqueness for
\begin{equation}\label{eqTransformee}
\bX^\obl (t)=\bX^\obl (0) +\int_0^t \obl^{-1} \bsigma(\obl \bX^\obl(s)) d\bW(s) +\int_0^t \obl^{-1} \bb(\obl \bX^\obl(s)) ds +\sum_{f\in\Const} \int_0^t~^t\obl \nabla f(\obl \bX^\obl(s)) dL_f(s)
\end{equation}
in the closure of the set $\obl^{-1} \D =\left\{ \by \in \R^d;~ \forall f\in\Const~f(\obl \by) > 0 \right\}$
with local times satisfying the condition $\disp L_f(\cdot)=\int_0^\cdot\un_{f(\obl \bX^\obl(s))=0}~dL_f(s)$.

The transformed coefficients $\bsigma^\obl=\obl^{-1} \bsigma(\obl~\cdot)$ and $\bb^\obl=\obl^{-1} \bb(\obl~\cdot)$ inherit the boundedness and Lipschitz continuity property from $\bsigma$ and $\bb$.
Remark \ref{remStabilite} provides the compatibility of the set of transformed constraints
$\{ f(\obl~\cdot);~f\in\Const \}$. Moreover, (\ref{eqTransformee}) is an equation with normal reflection because $\nabla(f(\obl ~\cdot))= {^t\obl} \nabla f(\obl~\cdot)$. Thus equation (\ref{eqTransformee}), and then equation (\ref{eqPrincipale}), have a unique strong solution.

Moreover, if $\bsigma=\obl$, $\bb=-\frac{1}{2} \obl {^t\obl} \nabla\Phi$ and 
$\bX(0) \sim \un_{\D}(\bx) e^{-\Phi(\bx)} d\bx$ for some $\C^2$-function $\Phi$ with bounded derivatives, 
then $\bX^\obl$ is the solution of equation (\ref{eqTransformee}) with $\bsigma^\obl=I_d$, 
$\bb^\obl=-\frac{1}{2} ^t\obl \nabla\Phi(\obl~\cdot)$ and initial distribution 
$\bX^\obl(0) \sim  \un_{\obl^{-1} \D}(\by) e^{-\Phi(\obl\by)} |det(\obl)| d\by$. 
Thanks to the reversibility result obtained for normal reflection, $\bX^\obl$ is time-reversible. This implies the time-reversibility of the solution of equation (\ref{eqPrincipale}).
\end{proof2}

%>>>>>>>>>>>>>>>>>>>>>>>>>>>>>>>>>>>>>>>>>>>>>>>>>>>>>>>>>>>>>>>>>>>>>>>>>>>>>>>>>
\section{An application : cluster of particles near the surface of an attractive planet}
\label{sect_examples}
%>>>>>>>>>>>>>>>>>>>>>>>>>>>>>>>>>>>>>>>>>>>>>>>>>>>>>>>>>>>>>>>>>>>>>>>>>>>>>>>>>
%>>>>>>>>>>>>>>>>>>>>>>>>>>>>>>>>>>>>>>>>>>>>>>>>>>>>>>>>>>>>>>>>>>>>>>>>>>>>>>>>>

As an application of the previous results, let us study the configuration of a large number of particles around a fixed sphere called planet.
As in \cite{FradonBrownianGlobules}, we consider spherical hard particles with random radii oscillating between a minimum and a maximum value. Each particle is driven by a Brownian motion, and undergoes the influence of the gravitational attraction generated by the planet. The motion is perturbated as the particles bump into each other and into the planet.
What we are interested in is the existence and uniqueness of such a dynamics, and the equilibrium distribution of the particles:
Do they tend to be scattered in a large area, or are they typically close together, all located beneath some altitude, with high probability~? That is, do they form some sort of "ocean" around the planet~?
Using the results of section \ref{sect_Skorohod_eq}, we prove in proposition \ref{ocean_globules} that the particles eventually tend to cluster at the surface of the planet when the temperature (represented by the diffusion coefficient) tends to zero. 

More precisely, the planet is the $d$-dimensional closed ball $B(0,R)$ centered at the origin with radius $R$.
The space around this planet contains a large number $n$ of particles, which are spheres with random positions and random radii. Each one will be represented by the position $x_i$ of its center in $\R^d$ and the value $\nx_i$ of its radius. Thus configurations are vectors $\bx=(x_1,\nx_1,\ldots,x_n,\nx_n)$ in $\R^{n(d+1)}$.  

To prevent negative radii, we enforce $\nx_i\in[\rm,\rp]$ for some fixed values $0<\rm<\rp$. 

The random oscillations of the positions of the particles are not on the same scale as the random oscillations of their radii.
We take this into account via the elasticity coefficient $\nsigma>0$ of their surface. 

We assume that the gravity field $\grav$ generated by the planet is isotropic: it only depends on the norm $|x|$ (more general cases would be mathematically tractable, but without physical meaning).
As usual (see e.g. \cite{BurdzyArchimede}) the gravitational attraction appears through a drift in the dynamics:
% to get around the problem that diffusions do not have pathes smooth enough for gravity to appear as a second derivative : the potential energy $-\grav'$ will directly influence the "speed" of the particle, rather than having its derivative (the gravity force) acting on the "acceleration" of the particle.
Function $\grav$ is an increasing function which is $\C^2$ on $]0;+\infty[$.
%We also assume that the gravitational attraction 
The drift decreases with the distance, but not too fast (or equivalently the potential energy $-\grav'$ increases not too fast) in the sense that $\grav'' \le 0$ and $~\liminf_{\rho\to +\infty}\rho\grav'(\rho) >0$.

An important example in dimension $d=3$ is $\grav(\rho)=C^{st}\ln(\rho)$ which gives the drift $-\grav'(\rho)=-\frac{C^{st}}{\rho}$ corresponding to the usual gravitational acceleration $-\grav''(\rho)=\frac{C^{st}}{\rho^2}$. 

\begin{figure}[!h]
  \centering
  \includegraphics[scale=0.9]{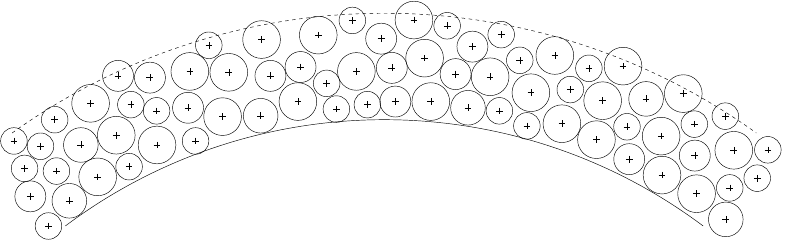}
  \caption{A configuration with an interface between high particle density and empty space.}\label{fig_ocean}
\end{figure}

At temperature $\theta>0$, the random motion of particles is modelized by the stochastic differential system:
$$
(\E_\theta)
\left\{
\begin{array}{l} 
  \text{for }i\in\{1,\ldots,n\}\\  \disp
  X_i(t) = X_i(0) +\theta W_i(t) -\int_0^t \grav'(|X_i(s)|) \frac{X_i}{|X_i|}(s)~ds \\ \disp
                                 \phantom{X_i(t) = X_i(0) +\theta W_i(t)}
                                 +\int_0^t \frac{X_i}{R+\nX_i}(s)~dL^R_i(s) 
                                 +\sum_{j=1}^n \int_0^t \frac{X_i-X_j}{\nX_i+\nX_j}(s)~dL_{ij}(s) \\ \disp
 \nX_i(t) = \nX_i(0) +\theta\nsigma \nW_i(t) -\nsigma^2 L^R_i(t) -L^+_i(t) +L^-_i(t) -\nsigma^2\sum_{j=1}^n L_{ij}(t)
\end{array}
\right.
$$
In this equation, vector $(X_i(\cdot),\nX_i(\cdot))_{1 \le i\le n}$ represents the positions and radii of the $n$ particles, the $W_i$'s are independent $\R^d$-valued Brownian motions and the $\nW_i$'s are independent one-dimensional Brownian motions, also independent from the $W_i$'s. The amplitude of the Brownian oscillation of the positions depends on temperature $\theta$, while the amplitude of the radii oscillation depends on both the temperature $\theta$ and the surface elasticity $\nsigma$.
Since $\grav'$ is positive, the drift of $X_i$ is always directed toward the origin, i.e. toward the planet as expected.
The local time $L^R_i$ represents the repulsion received by the $i^\text{th}$ particle when it collides with the planet (impulsion away from the origin, in the direction of the unit vector $\frac{X_i}{R+\nX_i}$), and the local times $L_{ij}$ represent the collisions between particles, which tend to move the involved particles away from each other: unit direction $\frac{X_i-X_j}{\nX_i+\nX_j}$.
Collisions between particles are symmetric ($L_{ij} \equiv L_{ji}$).These local times also appear in the dynamics of the radii, because particles are like bubbles, and collisions have a deflating effect on them (the radii decrease). The local times $L^+_i$ and $L^-_i$ are here to comply with the condition $\nx_i\in[\rm,\rp]$ and give a positive (resp. negative) impulsion to the radius if it becomes too small (resp. too large).
The impulsions are only given on the boundary of the set of allowed configurations, therefore the $L^R_i$'s, $L^+_i$'s, $L^-_i$'s and $L_{ij}$'s should satisfy:
$$
(\E'_\theta)
\left\{
\begin{array}{l} 
\text{for }i,j\in\{1,\ldots,n\}\\  \disp
L^R_i(t) = \int_0^t \un_{|X_i(s)|=R+\nX_i(s)} ~dL^R_i(s)~, \quad L^+_i(t) = \int_0^t \un_{\nX_i(s)=\rp} ~dL^+_i(s) \\ \disp
L^-_i(t) = \int_0^t \un_{\nX_i(s)=\rm} ~dL^-_i(s)~, \quad L_{ij}(t) = \int_0^t \un_{|X_i(s)-X_j(s)|=\nX_i(s)+\nX_j(s)} ~dL_{ij}(s)
\end{array}
\right.
$$ 
The set of constraints is therefore:
\begin{itemize}
\item $f_i^R(\bx)=|x_i|^2-(R+\nx_i)^2 >0$ for $1 \le i \le n$ \hfill(the particles do not intersect the planet);
\item $f_i^+(\bx)=\rp-\nx_i >0$ for $1 \le i \le n$ \hfill(radii are smaller than the maximum value); 
\item $f_i^-(\bx)=\nx_i-\rm >0$ for $1 \le i \le n$ \hfill(radii are larger than the minimum value); 
\item $f_{ij}(\bx)=|x_i-x_j|^2-(\nx_i+\nx_j)^2 >0$ for $i \neq j$ in $\{1,2,\ldots,n\}$ \hfill(particles do not overlap).
\end{itemize}
\begin{propal}\label{compatibilite_globules}~\\
%This collection 
$\disp \left\{ f_i^R, f_i^+, f_i^-;~1\le i \le n \right\}\cup\left\{ f_{ij};~1\le i < j \le n \right\}$ is a set of compatible constraints on $\R^{n(d+1)}$.
\end{propal}
Let $\D=\bigcap_{i=1}^n \left( (f_i^R)^{-1}(\R^*_+) \cap (f_i^+)^{-1}(\R^*_+) \cap (f_i^+)^{-1}(\R^*_+) 
                                     \cap \bigcap_{j \neq i} (f_{ij})^{-1}(\R^*_+) \right)$.
We want to study the effect of the attraction of the planet on a "typical" particle configurations, i.e. at equilibrium.             
\begin{propal}\label{solution_globules}~\\
If $\grav$ is an increasing $\C^2$-function on $]0;+\infty[$ satisfying $\grav'' \le 0$ and $~\disp \liminf_{\rho\to +\infty}\rho\grav'(\rho) >0$, then equation $(\E_\theta,\E'_\theta)$ has a unique strong solution, which is a $\overline{\D}$-valued process. \\
The measure $\un_{\D}(\bx) e^{-\frac{1}{\theta^2} \sum_{i=1}^n \grav(|x_i|)} d\bx$ is a time-reversible measure for the solution. For $\theta$ small enough, this measure is finite thus the solution admits a time-reversible probability measure:
$$
\mu_\theta(d\bx)=\frac{1}{\int_\D e^{-\frac{1}{\theta^2} \sum_{i=1}^n \grav(|y_i|)} d\by} 
                 \un_{\D}(\bx) e^{-\frac{1}{\theta^2} \sum_{i=1}^n \grav(|x_i|)} d\bx
$$
\end{propal}
Our aim, once the existence and uniqueness of the dynamics is proved, is to verify that, at low temperature, all particles cluster around the planet with high probability. 
In other words, there exists with high probability an interface between two regions around the planet: no particle over a certain altitude, and beneath this altitude a particle density so high that one cannot add any particle more (see figure \ref{fig_ocean}).
%This situation we characterize, in a way inspired by the nice ideas developped in \cite{BurdzyArchimede}, by the existence of some "surface" around the large sphere, with no particle above the surface and no empty space large enough to add one more particle below the surface:

\begin{propal}\label{ocean_globules}~\\
For each positive $\eps$, let $A_\eps$ be the set of configurations which do not pack into a minimal volume:
$$
A_\eps=\left\{ \bx\in\D;~\exists \by\in\D~ \exists k \le n \text{ s.t. } \forall i \neq k ~~ y_i=x_i \text{ and }
                                                                            |y_k|<|x_k|-\eps \right\}
$$
The probability that $A_\eps$ occurs at equilibrium tends to zero as the temperature tends to zero:
$$\disp \lim_{\theta\to 0} \mu_\theta(A_\eps)=0.$$
\end{propal}
%{\sc et prouver que $\mu_\theta(A_\eps) $ tend vers 1 si $\theta\to +\infty$ ??????????????????}
The end of the paper is devoted to the proofs of the three above propositions.

\begin{proof2}{\bf of proposition \ref{compatibilite_globules}}
%The constraints in $\disp \Const=\left\{ f_{ij}~;~1\le i < j \le n \right\}\cup \left\{ f_i^+, f_i^-, f_i^R~;~1\le i \le n \right\}$ are obviously $\C^2$. 

%The set $\D$ of possible particle configurations is a non-empty connected set because each possible configuration can be shifted along a continuous path in $\D$ to obtain a fixed configuration, say $((R+4\rp)e_1,\rm,(R+8\rp)e_1,\rm,(R+12\rp)e_1,\rm,\ldots,(R+4n\rp)e_1,\rm)$ where $e_1=(1,0,\cdots,0)$ denotes the first unit vector in $\R^d$:  pick the particle $(x_i,\nx_i)$ farthest from the origin, deflate it to $(x_i,\rm)$, then push it away from the origin and slide it on the sphere centered at $0$ with radius $|x_i|+4n^2\rp$ to the position $(|x_i|+4n^2\rp)e_1$, and repeat the procedure to slide the next particle farthest from the origin $(x_j,\nx_j)$ to the position $((|x_j|+4n(n-1)\rp)e_1,\rm)$, and so on, until all particles are on a straight line at distance at least $R+4n\rp$ from the origin; then slide down the particle labeled $1$ around the other particles to the position $(R+4\rp)e_1$ and do it again for each other particles.

The constraints in $\disp \Const=\left\{ f_{ij}~;~1\le i < j \le n \right\}\cup \left\{ f_i^+, f_i^-, f_i^R~;~1\le i \le n \right\}$ are $\C^\infty$ and the corresponding set $\D$ of possible configurations is obviously a non-empty connected set.
The first derivative of each constraint function is uniformly positive on its vanishing set because:
\begin{itemize}
\item $\nabla f_i^R(\bx)=2\left( 0,\ldots,0,x_i,-(R+\nx_i),0,\ldots,0 \right)$ \\
      if $f_i^R(\bx)=0$ i.e. $|x_i|=R+\nx_i$ then $|\nabla f_i^R(\bx)|=2\sqrt{2}(R+\nx_i) \ge 2\sqrt{2}(R+\rm) >0$~;
\item $\nabla f_i^+(\bx)=-\nabla f_i^-(\bx)=-\left( 0,\ldots,0,1,0,\ldots,0 \right)$ 
      \quad {\small ($(i(d_1+d_2+1)-1)$th coordinate)}~;
\item $\nabla f_{ij}(\bx)= 2\left( 0,\ldots,0,x_i-x_j,-(\nx_i+\nx_j),0,\ldots,0,x_j-x_i,-(\nx_i+\nx_j),0,\ldots,0 \right)$ \\
      if $f_{ij}(\bx)=0$ i.e. $|x_i-x_j|=\nx_i+\nx_j$ then $|\nabla f_{ij}(\bx)|=4(\nx_i+\nx_j) \ge 8\rm >0$.
\end{itemize}
%The constraint functions are linear or quadratic, their second derivatives are constant.
In order to check the condition $\disp \inf_{\bx\in\partial\D} d(0,\Conv(\bx)) >0$, in the form given in remark \ref{remEquivalence}, we have to find some positive $\beta_0$, and some non-vanishing vector $\bv$ depending on $\bx\in\partial\D$, such that
$$
\forall f\in\Const \text{ s.t. } f(\bx)=0 \quad \bv.\nabla f(\bx) \ge \beta_0 |\bv|~|\nabla f(\bx)|
$$
From an intuitive point of view, $\bv$ is the "shortest way to go back" into $\D$ from the point $\bx$ on the boundary of $\D$, i.e. the quickest way for colliding particles to go apart, for particles with maximum (resp. minimum) radius to become smaller (resp. larger) and for particles touching the planet to go away. $C^R$ will denote the indices of these globules: $C^R=\{i \text{ s.t. } |x_i|=R+\nx_i\}$.

Intuitively, the best way to separate colliding particles is to move them away from the center of gravity of the cluster, that is to give each center $x_i$ an impulsion in the direction $x_i-\frac{1}{\sharp C_i} \sum_{j\in C_i} x_j$ where $C_i\subset\{1,\ldots,n\}$ is the cluster of colliding particles around $x_i$ (i.e. $C_i$ is the set containing $i$ and all indices connected to $i$ in the graph constructed on the vertices $\{1,\ldots,n\}$ by the edges $j \sim j' \Longleftrightarrow |x_j-x_{j'}|=\nx_j+\nx_{j'}$). Similarly, the best way for particles touching the planet to go away is for each center $x_i$ with $i\in C^R$ to receive a small impulsion proportional to $x_i$ (this impulsion will also separate clusters of colliding particles). So a convenient $\bv$ should be:
$$
v_i = \left\{ \begin{array}{ll} 
              \disp x_i-\frac{1}{\sharp C_i} \sum_{j\in C_i} x_j & \text{if } C_i \cap C^R = \emptyset \\
                    x_i & \text{if } C_i \cap C^R \neq \emptyset 
              \end{array}
       \right. \quad\text{ and }\quad
\nv_i=\left\{ \begin{array}{ll} 
              \disp \phantom{-}\rm/2 & \text{if } \nx_i=\rm   \\
              \disp -\rm/2 & \text{if } \nx_i=\rp \\
              0 \text{ otherwise}
              \end{array}
       \right.
$$
Let us prove that the above vector $\bv$ satisfies the desired inequalities:
\begin{itemize}
\item if $|x_i|=R+\nx_i$ then $v_i=x_i$ hence:\quad
$\disp
\bv.\frac{\nabla f_i^R(\bx)}{|\nabla f_i^R(\bx)|} 
%=\frac{2}{2\sqrt{2}(R+\nx_i)}( x_i.v_i - (R+\nx_i)\nv_i )
=\frac{R+\nx_i}{\sqrt{2}}-\frac{\nv_i}{\sqrt{2}}
%\ge \frac{R+\rm}{\sqrt{2}}-\frac{\rm}{2\sqrt{2}}
\ge \frac{R}{\sqrt{2}}
$
\item if $\nx_i=\rp$ then $\disp \bv.\frac{\nabla f_i^+(\bx)}{|\nabla f_i^+(\bx)|}=-\nv_i=\frac{\rm}{2}$~; 
and if $\nx_i=\rm$ then $\disp \bv.\frac{\nabla f_i^-(\bx)}{|\nabla f_i^-(\bx)|}=\nv_i=\frac{\rm}{2} $~;
\item If $|x_i-x_j|=\nx_i+\nx_j$ then $C_i=C_j$ which implies $v_i-v_j=x_i-x_j$ thus:
$$
\bv.\frac{\nabla f_{ij}(\bx)}{|\nabla f_{ij}(\bx)|}
%=\frac{1}{4(\nx_i+\nx_j)} \left( (x_i-x_j).(v_i-v_j)-(\nx_i+\nx_j)(\nv_i+\nv_j) \right)
=\frac{\nx_i+\nx_j}{4} -\frac{\nv_i+\nv_j}{4}
\ge \frac{\rm}{4}
$$
\end{itemize}
So $\bv.\frac{\nabla f(\bx)}{|\nabla f(\bx)|}$ is bounded from below, uniformly in $\bx\in\partial\D$ and $f\in\Const$ such that $f(\bx)=0$. To complete the proof of proposition \ref{compatibilite_globules}, it only remains to find a uniform upper bound for $|\bv|$:
$$
|\bv|^2 = \sum_{i=1}^n |v_i|^2 +\nv_i^2
= \sum_{i;~C_i \cap C^R \neq \emptyset} |x_i|^2 
  +\sum_{i;~C_i \cap C^R =\emptyset} |\frac{1}{\sharp C_i} \sum_{j\in C_i} (x_i-x_j)|^2 +n\frac{\rm^2}{4}
$$
If $C_i$ is any cluster of colliding globules: 
$$
\left| \sum_{j\in C_i} (x_i-x_j) \right|^2 
\le \sharp C_i \sum_{j\in C_i} |x_i-x_j|^2 
\le \sharp C_i \sum_{k=0}^{\sharp C_i - 1} (2k\rp)^2 
= (2\rp)^2 (\sharp C_i)^2 \frac{(\sharp C_i - 1)(2\sharp C_i-1)}{6}
$$
Similarly, if $C_i$ is a cluster with at least one globule at distance $R+\nx_i$ of the origin: 
$$
\sum_{j\in C_i} |x_j|^2 
\le \sum_{k=0}^{\sharp C_i - 1} (R+\nx_i+2k\rp)^2 
\le 2\sharp C_i (R+\nx_i)^2 + 2(2\rp)^2 \frac{(\sharp C_i - 1)\sharp C_i(2\sharp C_i-1)}{6}
%\le 2n(R+\rp)^2 + \frac{4}{3} \rp^2 (n-1)n(2n-1)
$$
and the same upper bound holds for a sum over a union of such clusters. Consequently:
$$
|\bv|^2 \le 2n(R+\rp)^2 + \frac{4}{3} \rp^2 (n-1)n(2n-1)
            + \frac{2}{3}\rp^2 \sum_{i;~C_i \cap C^R = \emptyset} (\sharp C_i - 1)(2\sharp C_i-1) +n\frac{\rm^2}{4}
$$
%$$
%\le 2n(R+\rp)^2 + \frac{4}{3} \rp^2 (n-1)n(2n-1) + \frac{2}{3}\rp^2 n(n-1)(2n-1) +n\frac{\rm^2}{4}
%= 2n(R+\rp)^2 + 2\rp^2 (n-1)n(2n-1) + n\frac{\rm^2}{4}
%$$
Since the sum over $\{ i;~C_i \cap C^R = \emptyset \}$ is smaller than $n(n-1)(2n-1)$,
the norm of $\bv$ is uniformly bounded from above as a function of $\bx$. This completes the proof.
\end{proof2}

\begin{proof2}{\bf of proposition \ref{solution_globules}}\\
We use theorem \ref{thPrincipal} with the $n(d+1) \times n(d+1)$ diagonal matrix $\obl$ which has $n$ times the sequence $(\theta,\ldots,\theta,\theta\nsigma)$ as its main diagonal entries.
%$\obl$ is invertible for each positive temperature $\theta$.
Since the constraints are compatible on $\R^{n(d+1)}$, for any $\C^2$-function $\Phi$ on $\R^{n(d+1)}$ with bounded derivatives, the equation:
\begin{equation}\label{eqgeneric}
\bX(t)=\bX(0) + \obl \bW(t) -\frac{1}{2} \int_0^t \obl~^t\obl \nabla\Phi(\bX(s)) ds + \sum_{f\in\Const} \int_0^t \obl~^t\obl \nabla f(\bX(s)) dL_f(s) 
\end{equation}
has a unique strong solution in the closure of the set $\D$ defined by the constraints. Choosing 
$\Phi(\bx)=\sum_{i=1}^n \grav(|x_i|) / \theta^2$ hence $\nabla_{x_i}\Phi(\bx)=\frac{1}{\theta^2} \frac{x_i}{|x_i|} \grav'(|x_i|)$ and $\nabla_{\nx_i}\Phi(\bx)=0$, equation (\ref{eqgeneric}) becomes
{\small
$$
\!\!\!\!\!\!\!\!\left\{
\begin{array}{l} \disp
  X_i(t) = X_i(0) + \theta W_i(t) -\int_0^t \grav'(|X_i(s)|) \frac{X_i}{|X_i|}(s)ds +\int_0^t 2\theta^2 X_i(s)dL_{f_i^R}(s)
           +\sum_{j=1}^n \int_0^t 2\theta^2(X_i-X_j)(s)dL_{f_{ij}}(s) \\ \disp
 \nX_i(t) = \nX_i(0) + \theta\nsigma \nW_i(t) +\theta^2\nsigma^2\left(
                                                                -\int_0^t 2(R+\nX_i)(s)dL_{f_i^R}(s) 
                                                                -L_{f^+_i}(t) +L_{f^-_i}(t) 
                                                                -\sum_{j=1}^n \int_0^t 2(\nX_i+\nX_j)(s)dL_{f_{ij}}(s)
                                                                \right)
\end{array}
\right.
$$
}
Let us define $L_{ij}(\cdot) = 2\theta^2 \int_0^\cdot (\nX_i+\nX_j)(s)dL_{f_{ij}}(s)$, $L^+_i = \theta^2\nsigma^2L_{f^+_i}$, 
$L^-_i = \theta^2\nsigma^2L_{f^-_i}$ and $L^R_i(\cdot) = 2\theta^2 \int_0^\cdot (R+\nX_i)(s)dL_{f_i^R}(s)$.
The property $\disp L_f(\cdot) = \int_0^\cdot \un_{f(\bX(s))=0}~dL_f(s)$ 
%that each local time $L_f$ in (\ref{eqgeneric}) increases only on the level set $\{ f=0 \}$ 
implies that condition $(\E'_\theta)$ is satisfied for these new local times.
Then the solution of equation (\ref{eqgeneric}) is the solution of $(\E_\theta)$.
Thanks to theorem \ref{thReversibilite}, the measure $\un_{\D}(\bx) e^{-\Phi(\bx)} d\bx=\un_{\D}(\bx) e^{-\frac{1}{\theta^2} \sum_{i=1}^n \grav(|x_i|)} d\bx$ is a time-reversible measure for the solution.
To complete the proof, let us check that this measure can be renormalized as a probability measure for $\theta$ small enough.

Since $\ell:=\liminf_{\rho\to +\infty}\rho\grav'(\rho)>0$:
$\disp \quad \forall\eta>0 \quad \exists K >0 \quad \forall\rho>K \quad \grav'(\rho)\ge\frac{\ell-\eta}{\rho}$\\
For $\rho\ge K$ this integrates into $\disp \grav(\rho)\ge\grav(K)+(\ell-\eta)(\ln\rho-\ln K)$, so for any positive constant $c$:
$$
\int_{K}^{+\infty} e^{-c\grav(\rho)} \rho^{d-1} d\rho 
\le e^{-c\grav(K)} K^{c(\ell-\eta)} \int_{K}^{+\infty} \rho^{-c(\ell-\eta)+d-1} d\rho 
$$
For $c$ large enough to satisfy $-c(\ell-\eta)+d<0$ the above integral is finite, that is:
$$
\int_{\R^d\setminus B(0,R)} e^{-c\grav(|x|)} dx <+\infty
$$
This leads to the desired normalization constant, for $\theta$ small enough to satisfy $\frac{1}{\theta^2}\ge c$: 
$$
\begin{array}{rcl}\disp
\int_\D e^{-\frac{1}{\theta^2} \sum_{i=1}^n \grav(|x_i|)} d\bx
%&\le& \disp e^{-\frac{n}{\theta^2}\underline{\grav}} 
%            \int_\D e^{-\frac{1}{\theta^2}(\sum_{i=1}^n \grav(|x_i|)-n\underline{\grav})} d\bx \\ \disp
&\le& \disp e^{-\frac{n}{\theta^2}\underline{\grav}} 
            \int_{(\R^d\setminus B(0,R))^n} e^{-c(\sum_{i=1}^n \grav(|x_i|)-n\underline{\grav})} d\bx \\ \disp
&\le& \disp e^{c n\underline{\grav}-\frac{n}{\theta^2}\underline{\grav}} 
            \left( \int_{\R^{d-1}\setminus B(0,R)} e^{-c\grav(|x|)} dx \right)^n
            <+\infty      
\end{array}      
$$
where $\underline{\grav}=\min_{[R;+\infty[}\grav$ denotes the infimum on $[R;+\infty[$ of the smooth increasing function $\grav$.
%The existence of the normalization constant needed for the reversible measure to be a probability measure is proved, for $\theta<\sqrt{\ell/d}$.
\end{proof2}

\begin{proof2}{\bf of proposition \ref{ocean_globules}}\\
Let $\underline{\grav}_\D := \inf\{\sum_{i=1}^n \grav(|y_i|);~\by\in\D\}$. This infimum exists because $\grav$ is increasing on $]0;+\infty[$. We fix $\bx\in A_\eps$. There exists an allowed configuration $\by$ with all particles at the same position as in $\bx$ except one particle (say, the $k$'th) which satisfies $|y_k|<|x_k|-\eps$. Since $\grav'$ is a decreasing function:
$$
\sum_{i=1}^n \grav(|x_i|)=\sum_{i=1}^n \grav(|y_i|)+\int_{|y_k|}^{|x_k|} \grav'(\rho) d\rho
                         > \underline{\grav}_\D +(|x_k|-|y_k|)\grav'(|x_k|)
$$
Function $\grav'$ admits a limit at infinity:
\begin{itemize}
\item If this limit does not vanish, then
$\disp \sum_{i=1}^n \grav(|x_i|)> \underline{\grav}_\D +\eps\lim_{\rho\to +\infty} \grav'(\rho) >\underline{\grav}_\D$~;
\item If $\lim_{+\infty} \grav'=0$, the positivity of $\disp \ell:=\liminf_{\rho\to +\infty}\rho\grav'(\rho)$ implies the existence of a $K$ such that:
$$
\forall\rho\ge K \quad\quad \rho\grav'(\rho)\ge\frac{2\ell}{3} \quad\text{and}\quad (R+\rp)\grav'(\rho)\le\frac{\ell}{3} 
\quad\text{hence}\quad (\rho-R-\rp)\grav'(\rho)\ge\frac{\ell}{3}
$$
Without loss of generality, we can choose $K \ge R+2n\rp+n\eps$.
  \begin{itemize}
  \item If, among the $x_k$'s such that there exists $\by\in\D$ satisfying $|y_k|<|x_k|-\eps$ and $y_i=x_i$ for $i \neq k$, at least one has a norm smaller than $K$, then 
  $\disp \sum_{i=1}^n \grav(|x_i|) > \underline{\grav}_\D+\eps \grav'(K) > \underline{\grav}_\D$~;
  \item If no "$\eps$-pushable" $x_k$ has a norm smaller than $K$, then all particles in $\bx$ are at distance at least $K$ from the origin (because $K \ge R+2n\rp+n\eps$ and it is impossible to fill so large a volume with only $n$ particles without one of them to be "$\eps$-pushable"). When the $x_k$ which has the largest norm is shifted at distance $R+\rp$ from the origin, becoming $y_k$ and giving the configuration $\by$, we have
  $\disp
  \sum_{i=1}^n \grav(|x_i|) > \underline{\grav}_\D+(|x_k|-R-\rp)\grav'(|x_k|) \ge \underline{\grav}_\D+\frac{\ell}{3} >\underline{\grav}_\D
  $.
  \end{itemize}
\end{itemize}
So we obtain $\disp \sum_{i=1}^n \grav(|x_i|) > \underline{\grav}_\D+\eps'$ for all $\bx\in A_\eps$,
with a positive $\eps'$ equal to $\disp \eps\lim_{+\infty} \grav'$ if this limit does not vanish and to $\min(\eps\grav'(K),~\frac{\ell}{3})$ otherwise.

An immediate consequence is \quad
$\disp \mu_\theta(A_\eps) \le \mu_\theta(\{ \bx\in\D;~\sum_{i=1}^n \grav(|x_i|) > \underline{\grav}_\D+\eps' \})$.\\
The normalization constant of the probability measure $\mu_\theta$ is larger than:
$$
\int_\D \un_{\sum_{i=1}^n \grav(|x_i|) \le \underline{\grav}_\D+\eps'} 
                  ~e^{-\frac{1}{\theta^2} \sum_{i=1}^n \grav(|x_i|)} ~d\bx
\ge e^{-\frac{1}{\theta^2}(\underline{\grav}_\D+\eps')} 
    \int_\D \un_{\sum_{i=1}^n \grav(|x_i|) \le \underline{\grav}_\D+\eps'}~d\bx                  
$$
thus\quad 
$\disp
\mu_\theta(A_\eps) 
\le \frac{\int_\D \un_{\sum_{i=1}^n \grav(|x_i|) > \underline{\grav}_\D+\eps'} 
                  ~e^{-\frac{1}{\theta^2} (\sum_{i=1}^n \grav(|x_i|)-\underline{\grav}_\D-\eps')} ~d\bx}
         {\int_\D \un_{\sum_{i=1}^n \grav(|x_i|) \le \underline{\grav}_\D+\eps'} ~d\bx}                 
$.\\
The denominator does not depend on $\theta$. Dominated convergence theorem ensures that the numerator converges to zero when $\theta$ tends to $0$. So we obtain $\disp \lim_{\theta\to 0} \mu_\theta(A_\eps)=0$.
\end{proof2}


\begin{thebibliography}{99}

\bibitem{BurdzyArchimede}
K. Burdzy, Z. Chen and S. Pal,
{\em Archimedes' principle for Brownian liquid},
to appear in Annals of Applied Probability.

\bibitem{DupuisIshii}
P. Dupuis and H. Ishii,
{\em SDEs with oblique reflection on nonsmooth domains}, 
Annals of Probability {\bf 21} No. 1 (1993), 554-580.
{\em Correction : SDEs with oblique reflection on nonsmooth domains}, 
Annals of Probability {\bf 36} No. 5 (2008), 1992-1997.

\bibitem{FradonBrownianGlobules}
M. Fradon,
{\em Brownian dynamics of globules},
Elec. J. of Probability, {\bf 15} No. 6 (2010), 142-161.

\bibitem{FradonRoellyGlobulesInfinis}
M. Fradon and S. Roelly,
{\em Infinitely many {B}rownian globules with Brownian radii},
Stochastics and Dynamics {\bf 10} No. 4 (2010), 591-612.
\bibitem{IkedaWatanabe}
N. Ikeda and S. Watanabe,
Stochastic Differential Equations and Diffusion Processes.
North Holland - Kodansha (1981).

\bibitem{LionsSznitman}
P. L. Lions and A. S. Sznitman,
{\em Stochastic Differential Equations with Reflecting Boundary Conditions}, 
Com. Pure and Applied Mathematics {\bf 37} (1984), 511-537.

\bibitem{SaishoSolEDS}
Y. Saisho,
{\em Stochastic Differential Equations for multi-dimensional domain with reflecting boundary}, 
Probability Theory and Related Fields {\bf 74} (1987), 455-477.

\bibitem{SaishoTanakaBrownianBalls}      
Y. Saisho and H. Tanaka,
{\em Stochastic Differential Equations for Mutually Reflecting Brownian Balls},
Osaka J. Math. {\bf 23} (1986), 725-740.

\bibitem{SaishoTanemuraBrownianParticles}
Y. Saisho and H. Tanemura, 
{\em A {B}rownian ball interacting with infinitely many {B}rownian particles},
Tokyo J. Math. {\bf 16} (1993), 429-439.

\bibitem{SaishoTanakaSymmetry}      
Y. Saisho and H. Tanaka,
{\em On the symmetry of a Reflecting Brownian Motion Defined by Skorohod's Equation for a Multi-Dimensional domain},
Tokyo J. Math. {\bf 10} (1987), 419-435.

\bibitem{Skorokhod1et2}
A. V. Skorokhod
{\em Stochastic equations for diffusion processes in a bounded region 1, 2}, 
Theor. Veroyatnost. i Primenen. {\bf 6} (1961), 264-274;  {\bf 7} (1962), 3-23.

\bibitem{Tanaka}
H. Tanaka,
{\em Stochastic Differential Equations with Reflecting Boundary Conditions in Convex Regions},
Hiroshima Math. J. {\bf 9} (1979), 163-177.

\end{thebibliography}
\end{document}